\documentclass{amsart}

\usepackage{amsmath}
\usepackage{amsxtra}
\usepackage{amscd}
\usepackage{amsthm}
\usepackage{amsfonts}
\usepackage{amssymb}
\usepackage{mathtools}
\usepackage{eucal}
\usepackage[all]{xy}
\usepackage{graphicx}
\usepackage{mathrsfs}
\usepackage{color}
\usepackage{caption}
\usepackage{enumerate}
\usepackage{faktor}
\usepackage{bbm}
\usepackage[T1]{fontenc} 
\usepackage[utf8]{inputenc} 
\usepackage{tikz}
\usepackage{tikz-cd}
\usetikzlibrary{matrix,arrows,decorations.pathmorphing}

\usepackage{xcolor} 

\definecolor{WildStrawberry}{RGB}{238,41,103}
\definecolor{RedViolet}{RGB}{161,36,107}
\definecolor{Blue}{RGB}{45,47,146}

\usepackage{hyperref}
\hypersetup{
  colorlinks = true,
  urlcolor = WildStrawberry,
  citecolor = Blue,
  linkcolor = RedViolet
}

\theoremstyle{plain}
\newtheorem{theorem}{Theorem}
\newtheorem{lemma}[theorem]{Lemma}

\theoremstyle{definition}
\newtheorem{definition}[theorem]{Definition}

\theoremstyle{remark}

\numberwithin{equation}{section}

\newcommand{\Z}{\mathbb{Z}}

\newcommand{\calo}{\mathcal{O}}

\newcommand\restr[2]{{\left.\kern-\nulldelimiterspace#1\right|_{#2}}}

\makeatletter
\AtBeginDocument
 {
    \def\@thm#1#2#3{%
      \ifhmode
        \unskip\unskip\par
      \fi
      \normalfont
      \trivlist
      \let\thmheadnl\relax
      \let\thm@swap\@gobble
      \let\thm@indent\indent 
      \thm@headfont{\scshape}
      \thm@notefont{\fontseries\mddefault\upshape}%
      \thm@headpunct{.}
      \thm@headsep 5\p@ plus\p@ minus\p@\relax
      \thm@space@setup
      #1
      \@topsep \thm@preskip               
      \@topsepadd \thm@postskip           
      \def\dth@counter{#2}%
      \ifx\@empty\dth@counter
        \def\@tempa{%
          \@oparg{\@begintheorem{#3}{}}[]%
        }%
      \else
        \H@refstepcounter{#2}%
        \hyper@makecurrent{#2}%
        \let\Hy@dth@currentHref\@currentHref
        \AddToHookNext{para/begin}{\MakeLinkTarget*{\Hy@dth@currentHref}}%
        \def\@tempa{%
          \@oparg{\@begintheorem{#3}{\csname the#2\endcsname}}[]%
        }%
      \fi
      \@tempa
    }%
  \dth@everypar={%
    \@minipagefalse \global\@newlistfalse
    \@noparitemfalse
    \if@inlabel
      \global\@inlabelfalse
      \begingroup \setbox\z@\lastbox
       \ifvoid\z@ \kern-\itemindent \fi
      \endgroup
      \unhbox\@labels
    \fi
    \if@nobreak \@nobreakfalse \clubpenalty\@M
    \else \clubpenalty\@clubpenalty \everypar{}%
    \fi
  }%
}

\title[]{Small composite numbers in orbits of linear maps}

\author{Jose Reyes}
\address{Massachusetts Institute of Technology}
\email{aleu@mit.edu}

\date{\today}

\begin{document}

\begin{abstract}
  Generalized Cunningham chains are sets of the form $\{f^n(z)\}_{n\ge0}$ where all its elements are prime numbers and $f$ is a linear polynomial with integer coefficients. We generalize this definition further to include starting terms that are not prime, and we obtain the bound of $\ell(z)< z$ if $z$ is big enough, where $\ell(z)$ is the size of the generalized Cunningham chain. Unlike a direct generalization of previous results, which require $z$ to have a prime factor that does not divide the leading term of $f$, this result is only dependent on the size of $z$ and not on its prime factorization.
\end{abstract}

\maketitle

\setcounter{tocdepth}{2} 
\tableofcontents


\section{Introduction and basic definitions}\label{Intro}

Given an integer polynomial $P$, consider the set $P(\Z)=\{P(z):z\in \Z\}$. Which elements of it are prime numbers? Are there even any prime numbers in it at all?\bigskip

The Bouniakowsky conjecture encapsulates all these questions, as explained in \cite{BounConj}: given a polynomial $P$ of degree greater than 1, positive leading coefficient and with $$\gcd(P(1),P(2),\dots,P(z),\dots)=1,$$ it conjectures that there are infinitely many primes of the form $P(z)$ for $z\in\Z^+$.\bigskip

No polynomial is known that makes the conjecture true. For example, it is famously not known whether there are infinitely many primes of the form $z^2+1$ with $z\in\Z^+$.\bigskip

The issue becomes worse if we consider orbits of $P$, that is, sets of the form $\calo_P(z)=\{P^i(z):i\ge 1\}$ where $$P^i(z)=\underbrace{P(P(\dots(P}_{i\text{ times}}(z))\dots))$$ with $P\in\Z[x]$ and $z\in\Z$. Consider the specific case $P(z)=z^2-2z+2$ and an initial value of $z=3$. $\calo_P(3)$ is the set of Fermat numbers, which have a long history of study. Despite this, we do not know if there are an infinite number of Fermat primes, and more importantly, we do not even know if there are an infinite number of composite Fermat numbers (see \cite{SumFermat} for a list of known factorizations of Fermat numbers), which is more surprising given the abundance of composite numbers in comparison to prime numbers.\bigskip

In this paper we will focus on the simple case of linear polynomials. A specific case that is of significant interest is the polynomial $P(z)=2z+1$. Prime numbers $p$ such that $2p+1$ is also prime are called \textit{Sophie Germain primes}, and lists $$\{p,2p+1,2(2p+1)+1,\dots,2^n(p+1)-1\}$$ of orbits using the polynomial $2z+1$ starting at a prime number such that all terms of the list are prime are called \textit{Cunningham chains}. Löh conducted an exhaustive search for Cunningham chains with starting term less than $2^{50}$ in \cite[Section 2. Computations]{CunChain}. A small example of a Cunningham chain is $\{41,83,167\}$.\bigskip

A \textit{generalized Cunningham chain} is a list $$\{p,f(p),f^2(p),\dots,f^{\lambda(p)-1}(p)\}$$ where all its elements are prime and $f(z)=az+b$ with $a,b$ relatively prime (note that this is necessary in order for $f(p)$ to be prime).\bigskip

A \textit{complete generalized Cunningham chain} (using the same $f$) is a list $$\{p,f(p),f^2(p),\dots,f^{\lambda(p)-1}(p)\}$$ that is a generalized Cunningham chain, has $f^{\lambda(p)}(p)$ composite and has $f^{-1}(p)$ composite or not an integer (so it cannot be extended forwards or backwards). We say this Cunningham chain has \textit{length} $\lambda(p)$.\bigskip

For $f(z)=2z+1$, one can prove that $\lambda(p)\le p-1$ with Fermat's little theorem and Kanado gives a better bound of $\displaystyle\lambda(p)<\frac{p}{2}$ in \cite[p.2]{FibChain}. We will repeat these same techniques in Sections \ref{MehSection} and \ref{Yaysection} and generalize them to all linear $f$ with positive coefficients. We will also expand the definition of Cunningham chains in order to cover a wider breadth of chains of primes.\bigskip

\begin{definition}\label{RootDef}
    Given integers $z,a,b$ with $a,b$ relatively prime, a \textit{rooted Cunningham chain} is a list $$\{f(z),f^2(z),\dots,f^{\ell(z)}(z)\}$$ such that all its elements are prime and $f^{\ell(z)+1}(z)$ is composite (here, $f(z)=az+b$). We say $\ell(z)$ is its length and that it has root $z$.
\end{definition}\bigskip

Notably, in this definition $z$ is excluded from the list itself, so it is not necessary for it to be prime. Definition 1 generalizes all non-generalized complete Cunningham chains (that is, complete Cunningham chains that use $f(z)=2z+1$) except $\{2,5,11,23,47\}$ because they are rooted Cunningham chains with root $\displaystyle\frac{p-1}{2}$, even if $\displaystyle \frac{p-1}{2}$ is not prime (since 2 is the only even prime). It is also important to mention that this definition of length is one less than if the root $z$ of the rooted Cunningham chain is included into the list (if $z$ is prime), which was our previous definition of length.\bigskip

The goal of this paper is to prove the following two elementary bounds on the length of a rooted Cunningham chain that generalize the bounds for Cunningham chains.\bigskip

\begin{theorem}\label{Meh}
    For coprime positive integers $a,b$ with $a>1$, define $f\colon \Z^+\rightarrow\Z^+$ as $f(z)=az+b$. For all $z\in \Z_{>1}$ coprime with $a$, the rooted Cunningham chain $$\{f(z),f^{2}(z),\dots,f^{\ell(z)}(z)\}$$ has length $\ell(z)<z$.
\end{theorem}\bigskip

\begin{theorem}\label{Yay}
    For coprime positive integers $a,b$ with $a>1$, define $f\colon \Z^+\rightarrow\Z^+$ as $f(z)=az+b$. There exists $M\in\Z^+$ such that for all $z\in\Z^+$ with $z>M$, the rooted Cunningham chain $$\{f(z),f^{2}(z),\dots,f^{\ell(z)}(z)\}$$ has length $\ell(z)<z$.
\end{theorem}\bigskip

The case $a=1$ transforms the Cunningham chain into an arithmetic progression. Results on the distribution of primes in an arithmetic progression are covered elsewhere in greater detail, such as in Dirichlet's theorem on arithmetic progressions or Green and Tao's paper \cite{MR2415379} which directly addresses lengths.\bigskip

Theorems \ref{Meh} and \ref{Yay} give the same bound with different conditions on $z$. Theorem \ref{Meh} is a direct generalization of the result on non-generalized Cunningham chains. However, it requires $z$ to be relatively prime to $a$, so there is an infinite subset of $\Z^+$ that is not addressed by Theorem \ref{Meh}. Theorem \ref{Yay} requires additional machinery but uses the same essential ideas as Theorem \ref{Meh}, and it provides the desired bound for all but finite $z\in \Z^+$. In both cases we have restricted the linear maps to $a,b$ positive. In the case of Theorem \ref{Meh}, this gets rid of some pesky cases such as the one where $z$ remains fixed by $f$, but the condition on the coefficients being positive can be relaxed and Theorem \ref{Meh} remains true. For Theorem \ref{Yay}, setting $a,b$ to be positive avoids the same issues, but if we relax similar conditions the bound becomes $\ell(z)<z+C$ for some constant $C>0$ dependent on $a$ and $b$, so we omit this case.\bigskip

In Section \ref{MehSection}, we will provide a proof of Theorem \ref{Meh} using Fermat's little theorem and some casework; in Section \ref{Yaysection} we will provide a proof of Theorem \ref{Yay} using two additional lemmas, which make use of a special sequence to generate new prime divisors.\bigskip

\section{Proof of Theorem \ref{Meh}}\label{MehSection}

Firstly, note that $f(z)=az+b>z+b>z$ so $f(z)>z$ always. We can also calculate $f^n(z)$ explicitly:

\begin{equation*}\begin{split}
    f^n(z) & =a^nz+b+ab+a^2b+\cdots+a^{n-1}b\\
    & =a^nz+b\cdot\frac{a^n-1}{a-1}
\end{split}
\end{equation*}\bigskip

\textit{Proof of Theorem \ref{Meh}:} If $z$ and $b$ share a common factor $d>1$, then $d$ divides both $a^nz$ and $\displaystyle b\cdot \frac{a^n-1}{a-1}$ while being smaller than their sum, so $f^n(z)$ is composite for all $n$. We may then assume that $z$ and $b$ are coprime for the rest of the proof.\bigskip

Since $z>1$, we may take a prime divisor $p$ of $z$. Because $z$ is coprime with $a$, $p$ does not divide $a$. We proceed with two cases, considering whether $p$ divides $a-1$ or not.\bigskip

\textit{Case 1:} $p$ divides $a-1$.

\begin{equation*}\begin{split}
    f^p(z) & =a^pz+b+ab+a^2b+\cdots+a^{p-1}b\\
    & \equiv a^p\cdot 0+b+1b+1^2b+\cdots+1^{p-1}b \text{ (mod $p$)}\\
    & \equiv pb \text{ (mod $p$)} \\
    & \equiv 0 \text{ (mod $p$)}
\end{split}
\end{equation*}\bigskip

Therefore, $p$ divides $f^p(z)$. Because $p$ is a divisor of $z$ and $f(z)>z$, $p\le z<f^p(z)$, so $f^p(z)$ is composite and the length $\ell(z)$ obeys $\ell(z)<p\le z$, as desired.\bigskip

\textit{Case 2:} $p$ does not divide $a-1$.\bigskip

Using Fermat's little theorem, $a^{p-1}-1$ is divisible by $p$. Because $a-1$ and $p$ are coprime and both divide $a^{p-1}-1$, their product divides $a^{p-1}-1$ as well, which is equivalent to $\displaystyle\frac{a^{p-1}-1}{a-1}$ being divisible by $p$. Then 

\begin{equation*}\begin{split}
    f^{p-1}(z) & =a^{p-1}z+b\cdot\frac{a^{p-1}-1}{a-1}\\
    & \equiv a^{p-1}\cdot 0+b\cdot 0 \text{ (mod $p$)}\\
    & \equiv 0 \text{ (mod $p$)} 
\end{split}
\end{equation*}\bigskip

Therefore, $p$ divides $f^{p-1}(z)$. Because $p$ is a divisor of $z$ and $f(z)>z$, $p\le z<f^{p-1}(z)$, so $f^{p-1}(z)$ is composite and the length $\ell(z)$ obeys $\ell(z)<p-1<z$, as desired.\qed\bigskip

The main idea of this proof of Theorem \ref{Meh} is that the simplest way to obtain that the expression $$f^n(z)=a^nz+b\cdot\frac{a^n-1}{a-1} $$ is composite is to find a common factor of $a^nz$ and $\displaystyle b\cdot\frac{a^n-1}{a-1}$. One may also notice that the only condition on the prime factor $p$ of $z$ was that it did not divide $a$, so we may relax the condition on $z$ in Theorem \ref{Meh} to $z$ having a prime factor that does not divide $a$ (instead of being relatively prime to $a$, which is stronger).\bigskip

The set of $z$'s that share all their prime factors with $a$ is always infinite: pick any prime factor $q$ of $a$ and $\{q,q^2,\dots,q^n,\dots\}$ is an infinite subset of $z$'s that have this property. This means that, while Theorem \ref{Meh} works on any $z$ relatively prime to $a$ regardless of size, it says nothing about an infinite subset of positive integers. On the other hand, Theorem \ref{Yay} works on all but a finite number of positive integers.\bigskip

One way of exploiting Theorem \ref{Meh} to eliminate the condition of $z$ being coprime with $a$ is to consider $f(z)$ instead of $z$, since $f(z)=az+b$ is always coprime with $a$ for $z\in\Z^+$. If we consider the rooted Cunningham chain $\{f^2(z),f^3(z),\dots,f^{\ell(f(z))} (z)\}$ with root $f(z)$, we can apply Theorem \ref{Meh} directly. This tells us $\ell(f(z))<f(z)$, which in turn implies $$\ell(z)\le\ell(f(z))+1<f(z)+1=az+b+1$$ by appending $f(z)$ to the beginning of the list $\{f^2(z),f^3(z),\dots,f^{\ell(f(z))} (z)\}$. We conclude that $\ell(z)<az+b+1$ in all cases, which is already a linear bound. Unfortunately, it depends on the size of both $a$ and $b$.\bigskip

Let us look at an example of how this proof of Theorem \ref{Meh} and subsequent discussion works. Take $f(z)=2z+3$ and consider the rooted Cunningham chain $\{5,13,29,61\}$ with root 1. Because Theorem \ref{Meh} requires the root of the Cunningham chain to have at least one prime factor, we cannot apply it directly. However, as we discussed in the previous paragraph, if we consider $f(1)=5$ we'll be able to leverage this prime factor to get a composite number further down the line, and because 5 is coprime with $a-1=1$, we know that this will be maximally at the term $f^4(5)$. Indeed, $f^4(5)=f(61)=125=5^3$.\bigskip

How can we prove Theorem \ref{Yay} if we require our prime factors to divide $z$?
We have been implicitly using the fact that linear maps are bijective modulo some prime $p$ by applying Fermat's Little Theorem. It is not guaranteed that repeated applications of $f$ will eventually return to 0 (mod $p$) if we do not \textit{start} with a term that is 0 (mod $p$). But it turns out that we will be able to do it by considering a specific sequence, as we will explain in the next section.

\section{Proof of Theorem \ref{Yay}}\label{Yaysection}

Our answer lies in the sequence $(s_n)_{n\ge 1}$:

\begin{definition}\label{s seq}
Define the sequence $(s_n)_{n\ge 1}$ explicitly as
\begin{equation*}
\begin{split}
    s_1 & =z-b\\
    s_2 & =z-b-ab\\
        & \;\;\vdots\\
    s_n &=z-b-ab-\cdots-a^{n-1}b\\
        &=z-b\cdot\frac{a^n-1}{a-1}\\
        &\;\;\vdots  
\end{split}
\end{equation*}

or inductively as $s_1=z-b$, $s_{n+1}=s_n-a^nb$.
\end{definition}

Consider a prime divisor $p$ of $s_n=z-b-ab-\cdots-a^{n-1}b$. We have

\begin{equation*}\begin{split}
    f^{m}(z) & =a^{m}z+b+ab+\cdots+a^{m-1}b\\
    & \equiv a^{m}\cdot (b+ab+\cdots+a^{n-1}b)+b+ab+\cdots+a^{m-1}b \text{ (mod $p$)}\\
    & \equiv b\cdot(1+a+a^2+\cdots+a^{m-1}+a^m+\cdots+a^{m+n-1}) \text{ (mod $p$)} \\
    & \equiv b\cdot\left(\frac{a^{m+n}-1}{a-1}\right) \text{ (mod $p$)}
\end{split}
\end{equation*}\bigskip

This is a similar expression to what we had before in the proof of Theorem \ref{Meh}. In fact, it is exactly the same as what we had when $p\mid z$, with the exponent on $a$ changed from $m$ to $n+m$. This is no coincidence: a nice feature of linear maps is that they have an explicit inverse that is also a linear map, that is, $$f^{-1}(z)=\frac{z-b}{a}.$$

We have $$f^{-n}(z)=\frac{z-b-ab-\cdots-a^{n-1}b}{a^n}=\frac{s_n}{a^n},$$ so we may suspect that we should somehow be able to include the numerator of this expression in our previous framework in order to obtain other composite numbers, and our suspicions would be correct.\bigskip

Before proving Theorem \ref{Yay}, we will examine the prime factorizations for elements of the sequence $(s_n)_{n\ge 1}$ (\ref{s seq}) closely with two crucial lemmas. \bigskip

\begin{lemma}\label{MainLemma}
Let $a$ be a positive integer with $k$ distinct prime factors.\\ If $z>b+ab+a^2b+\cdots+a^{k+1}b$, then there exists a prime $p$ that divides $s_i$ for some index $i$ with $1\le i\le k+1$ but does not divide $a$.
\end{lemma}\bigskip

In order to prove Lemma \ref{MainLemma}, we need a better description of the prime factorizations of the elements of the sequence $(s_n)_{n\ge 1}$ (\ref{s seq}). We will get this description from Lemma \ref{ValLemma}. Let $\nu_p(x)$ be the $p$-adic valuation of $x$ (the exponent of $p$ in the prime factorization of $x$).\bigskip

\begin{lemma}\label{ValLemma}
    If $p$ is a prime factor of $a$, then for each $n\in \Z^+$ we will have two possibilities: $$\begin{cases}
        \nu_p(s_{n+1})=n\nu_p(a) & \text{ if $\nu_p(s_n)$}> n\nu_p(a)\\
        \nu_p(s_{n+1})\ge \nu_p(s_n) & \text{ otherwise }
        
    \end{cases}$$
    Additionally, if $\nu_p(s_{n+1})=n\nu_p(a)$, then for every $n_1>n$ we will have $$\nu_p(s_{n_1})=n\nu_p(a)$$
\end{lemma}\bigskip

\textit{Proof of Lemma \ref{ValLemma}:} If $\nu_p(s_n)>n\nu_p(a)$, then

\begin{equation*}
\begin{split}
    \nu_p(s_{n+1}) &=\nu_p(s_n-a^nb)\\
    &=\min\{\nu_p(s_n),\nu_p(a^nb)\}\\ &=\min\{\nu_p(s_n),n\nu_p(a)\}\\
    &=n\nu_p(a)
\end{split}
\end{equation*}\bigskip

\noindent since $\gcd(a,b)=1$ means $\nu_p(b)=0$ and $\nu_p(s_n)>\nu_p(a^nb)=n\nu_p(a)$. We may prove by induction the claim that for every $n_1>n$ we will have $\nu_p(s_{n_1})=n\nu_p(a)$; the above is the base case $n_1=n+1$. For the inductive step, if $\nu_p(s_{n_1})=n\nu_p(a)$ for $n_1>n$, then

\begin{equation*}
\begin{split}
    \nu_p(s_{n_1+1}) &=\nu_p(s_{n_1}-a^{n_1}b)\\
    &=\min\{\nu_p(s_{n_1}),\nu_p(a^{n_1}b)\}\\ &=\min\{n\nu_p(a),n_1\nu_p(a)\}\\
    &=n\nu_p(a)
\end{split}
\end{equation*}\bigskip

\noindent since $\nu_p(s_{n_1})<\nu_p(a^{n_1}b)=n_1\nu_p(a)$, as desired.\bigskip

Otherwise, we'll have $\nu_p(s_n)\le n\nu_p(a)$, so $p^{\nu_p(s_n)}$ divides both $s_n$ and $a^nb$. Therefore, it divides their difference, $s_{n+1}$, which yields $\nu_p(s_{n+1})\ge \nu_p(s_n)$. \qed\bigskip

Lemma \ref{ValLemma} gives us a characterization of indices $n$ for which the $p$-adic valuation of $s_n$ breaks out of a weakly ascending pattern, in relation to a specific prime $p$. In order to formalize this we will use the following definition.\bigskip

\begin{definition}\label{nstable}
    For a prime $p$ that divides $a$, we say that $p$ is $n$\textit{-stable} if $$\nu_p(s_{n+1})=n\nu_p(a).$$
\end{definition}\bigskip

 Each prime that divides $a$ may only be $n$-stable for at most one $n$, because of the last sentence of Lemma \ref{ValLemma}. Lemma \ref{ValLemma} gives us all the necessary tools to prove Lemma \ref{MainLemma}.\bigskip

\textit{Proof of Lemma \ref{MainLemma}:} Assume the contrary. Then the prime factors of $s_i$ can only be some of the $k$ primes that divide $a$. By Lemma \ref{ValLemma} and the fact that each prime is $n$-stable for at most one $n$, there is at least one $1\le i\le k+1$ such that none of the $k$ distinct primes that divide $a$ are $i$-stable. For this specific $i$, we'll have that $\nu_p(s_{i+1})\ge\nu_p(s_i)$ for all $p$ that divide $a$.\bigskip

By our initial assumption, these are actually all the primes that are in the prime factorizations of $s_i$ and $s_{i+1}$, so this means that $s_i$ divides $s_{i+1}$. Therefore, $$s_i\mid (s_{i+1}-s_i)=a^ib$$ $$\Rightarrow (z-b-ab-\cdots-a^{i-1}b)\mid a^ib$$ $$\Rightarrow z-b-ab-\cdots-a^{i-1}b\le a^ib$$ $$\Rightarrow z\le b+ab+\cdots+a^ib\le b+ab+\cdots a^{k+1}b$$ which is a contradiction with our assumption on the size of $z$. \qed\bigskip

From Lemma \ref{MainLemma} we can essentially repeat the same proof we had for Theorem \ref{Meh} to prove Theorem \ref{Yay}, but with a prime divisor of this $s_i$ that does not divide $a$ instead of a prime divisor of $z$.\bigskip

\textit{Proof of Theorem \ref{Yay}:} If $z$ and $b$ share a common factor $f(z)$ is composite, so assume that $\gcd (z,b)=1$. Recall $a$ has $k$ distinct prime factors, and set $$M=b+ab+\cdots+a^{k+1}b.$$ If $z>M$, then by Lemma \ref{MainLemma} there exists a prime $p$ that divides $s_i$ (for some term of $(s_n)_{n\ge 1}$ from Definition \ref{s seq}) for some index $i$ with  $1\le i\le k+1$ but does not divide $a$.\bigskip

Because $z\equiv b+ab+\cdots+a^{i-1}b \text{ (mod $p$)}$, we'll have 

\begin{equation*}
\begin{split}
    f^n(z) &=a^nz+b+ab+\cdots+a^{n-1}b\\
    & \equiv a^n(b+ab+\cdots+a^{i-1}b)+b+ab+\cdots+a^{n-1}b \text{ (mod $p$)}\\
    & \equiv b\cdot (1+a+a^2+\cdots+a^{n+i-1})\text{ (mod $p$)}\\
    & \equiv b\cdot \left(\frac{a^{n+i}-1}{a-1}\right) \text{ (mod $p$)}
\end{split}
\end{equation*}\bigskip

We again proceed with two cases, considering whether $p$ divides $a-1$ or not.\bigskip

\textit{Case 1:} $p$ divides $a-1$.\bigskip

Let $r\in \{1,2,\dots,p\}$ be the residue of $-i\text{ (mod $p$)}$. We'll have 

\begin{equation*}\begin{split}
    f^r(z) & \equiv b\cdot (1+a+a^2+\cdots+a^{r+i-1}) \text{ (mod $p$)}\\
    & \equiv b\cdot (1^0+1^1+1^2+\cdots+1^{r+i-1}) \text{ (mod $p$)}\\
    & \equiv b\cdot(r+i) \text{ (mod $p$)} \\
    & \equiv 0 \text{ (mod $p$)}
\end{split}
\end{equation*}\bigskip

So $p$ divides $f^r(z)$. Because it is a divisor of $s_i$ and $f(z)>z$, $p\le s_i<z<f^r(z)$, so $f^r(z)$ is composite and the length $\ell(z)$ obeys $\ell(z)<r\le p<z$ as desired.\bigskip

\textit{Case 2:} $p$ does not divide $a-1$.\bigskip

Using Fermat's little theorem, $a^{(p-1)m}-1$ is divisible by $p$ for all $k\in\Z^+$. Because $a-1$ and $p$ are coprime and both divide $a^{(p-1)m}-1$, their product divides $a^{(p-1)m}-1$ as well, which is equivalent to $\displaystyle\frac{a^{(p-1)m}-1}{a-1}$ being divisible by $p$. \bigskip

Let $r\in\{1,2,\dots,p-1\}$ be the residue of $-i \text{ (mod $p-1$)}$. Then 

\begin{equation*}\begin{split}
    f^{r}(z) & \equiv b\cdot\frac{a^{r+i}-1}{a-1}\\
    & \equiv b\cdot 0 \text{ (mod $p$)}\\
    & \equiv 0 \text{ (mod $p$)} 
\end{split}
\end{equation*}

\noindent since $r+i$ is a multiple of $p-1$.\bigskip

So $p$ divides $f^{r}(z)$. Because $p$ is a divisor of $s_i$ and $f(z)>z$, $p\le s_i<z<f^{r}(z)$, so $f^{r}(z)$ is composite and the length $\ell(z)$ obeys $\ell(z)<r\le p-1<z$ as desired.
\qed\bigskip

Let us look at a final example that gives a tighter bound than Theorem \ref{Yay} by carefully considering $(s_n)_{n\ge 1}$.\bigskip

Let $f(z)=2z+3$ and consider the rooted Cunningham chain with root 32. We cannot apply Theorem \ref{Meh} because $32=2^5$. However, we can apply Theorem \ref{Yay} because $$32>M=b+ab+a^2b+\cdots+a^{k+1}b=3+2\cdot 3+2^2\cdot 3=21.$$ In this case, Theorem \ref{Yay} tells us $\ell(32)<32$, so the rooted Cunningham chain may not exceed 31 terms. Unfortunately, this is much bigger than the actual length of the chain, so let us look at $(s_n)_{n\ge 1}$ closely to see if we can obtain a better bound without directly computing the chain.\bigskip

The positive terms of $(s_n)_{n\ge 1}$ are 
\begin{equation*}
\begin{split}
    s_1 & =32-3 =29\\
    s_2 & =29-2\cdot 3 =23\\
    s_3 & =23-2^2\cdot 3 =11
\end{split}
\end{equation*}

\noindent since $s_4=11-2^3\cdot 3=-13$. Additionally, 11 is the term with least absolute value.\bigskip

Because $11$ does not divide $2-1=1$, Case 2 of the final part of the proof of Theorem \ref{Yay} applies. The residue $r$ of $-3\text{ (mod $10$)}$ is $7$, so 11 divides $f^7(32)=4477$. This means $\ell(32)\le 6$. \bigskip

This bound is actually tight: the full rooted Cunningham chain is $$\{67,137,277,557,1117,2237\}$$ with length 6.

\section{Acknowledgements}
I would like to thank Robin Zhang for leading the Spring 2025 Number Theory Seminar, from which the reason to work on this topic originated, and again to Robin Zhang, Lucy Epstein, Ray Wang and Ezra Guerrero for reviewing and making many (many) helpful annotations and suggestions throughout all steps of writing this paper. The content would be much less readable without them.


\bibliographystyle{alpha}
\bibliography{bibliography}

\end{document}